\documentstyle{article}

\newtheorem{thm}{Theorem}[section]
\newtheorem{cor}[thm]{Corollary}

\newtheorem{prop}[thm]{Proposition}
\newtheorem{defn}[thm]{Definition}
\newtheorem{rem}[thm]{Remark}

\begin{document}

\baselineskip = 15 pt
\parskip 4pt

\newtheorem{defin}{Definition}
\newcommand{\Al}{${\cal A}_{\hbar,\eta}(\widehat{sl}_2)~$}
\newcommand{\Alo}{${\cal A}_{\hbar,0}(\widehat{g})~$}
\newcommand{\Apq}{${\cal A}_{q,p}(\widehat{sl}_2)~$}
\newcommand{\Apqn}{${\cal A}_{q,p}(\widehat{sl}_N)~$}
\newcommand{\Apqg}{${\cal A}_{q,p}(\widehat{g})~$}
\newcommand{\Eg}{${\cal E}_{q,p}(\widehat{g})~$}
\newcommand{\Egc}[2]{${\cal E}_{q^{#1},p}(\widehat{g})_{#2}~$}
\newcommand{\Eegc}[2]{{\cal E}_{q^{#1},p}(\widehat{g})_{#2}}
\newcommand{\Alg}{${\cal A}_{\hbar,\eta}(\widehat{g})~$}
\newcommand{\Algi}[2]{{\cal A}_{\hbar,#1}(\widehat{g})_{#2}}
\newcommand{\Alp}[2]{{\cal A}_{#1,p}(\widehat{sl_2})_{#2}}
\newcommand{\Algp}[2]{{\cal A}_{#1,p}(\widehat{g})_{#2}}
\newcommand{\AlN}{${\cal A}_{\hbar,\eta}(\widehat{sl}_N)~$}
\newcommand{\phq}[4]{\phi^{(#1)}\left(\frac{#2}{#3}{#4}\right)}
\newcommand{\Phq}[4]{\Phi^{(#1)}_{ij}\left(\frac{#2}{#3}{#4}\right)}
\newcommand{\Psq}[4]{\Psi^{(#1)}\left(\frac{#2}{#3}{#4}\right)}

\rightline{math.QA/9801062}
\rightline{Revised April 5, 1998}
\vspace{2cm}

\centerline{\LARGE \bf Infinite Hopf Family of Elliptic Algebras and Bosonization}

\vspace{0.5cm}

\centerline{Bo-Yu Hou$^{b}$\footnote{E-mail: byhou@nwu.edu.cn}
\hspace{0.5cm} Liu Zhao$^{ab}$\footnote{E-mail: lzhao@nwu.edu.cn}
\hspace{0.5cm} Xiang-Mao Ding$^{ac}$\footnote{E-mail: xmding@itp.ac.cn}}

\centerline{${}^a$CCAST (World Lab), Academia Sinica,
P O Box 8730, Beijing 100080, P R China}
\centerline{${}^b$Institute of Modern Physics, Northwest University, Xian 710069,
China}
\centerline{${}^c$Institute of Theoretical Physics,
Academy of China, Beijing 100080, China}




\begin{abstract}
Elliptic current algebras ${\cal E}_{q,p}(\widehat{g})$
for arbitrary simply laced finite dimensional
Lie algebra $g$ are defined and their co-algebraic structures are studied. It is shown
that under the Drinfeld like comultiplications, the algebra ${\cal E}_{q,p}(\widehat{g})$
is not co-closed for any $g$. However putting the algebras ${\cal E}_{q,p}(\widehat{g})$
with different deformation parameters together, we can establish a structure of
infinite Hopf family of algebras. The level 1 bosonic realization for the algebra
${\cal E}_{q,p}(\widehat{g})$ is also established.
\end{abstract}

\section{Introduction}

In this paper we continue our recent study on infinite Hopf family of
algebras  and obtain new example of such families--
infinite Hopf family of elliptic algebras.

The concept of infinite Hopf family of algebras was first introduced our earlier paper
\cite{f} in which the algebras \Alg are proposed and their co-algebraic structures
are specified. In contrast to the standard Hopf structures for the quantum affine algebras
and Yangian doubles, the algebras \Alg, including their most degenerated case
\Al \cite{KLP}, are not evidently co-closed, and their co-algebraic structures
are formulated in terms of some generalized Hopf structure, examples are
the Hopf family of algebras of \cite{KLP} and the infinite Hopf family of algebras
of our paper \cite{f}.

The algebras \Alg appeared in \cite{f} are very unusual. For
$g=sl_2$, such algebra was proposed as the scaling limit
of the elliptic algebra \Apq and thus inherits {\em two}
deformation parameters $\hbar$ and $\eta$. The first parameter
$\hbar$ can be viewed as a ``quantization parameter'', because
in the limit $\hbar \rightarrow 0$, the algebra \Alg
would become a classical algebra. The second parameter $\eta$
should be viewed as a ``family deformation parameter'',
because the set of algebras \Alg with different $\eta$
form the first known nontrivial example of infinite Hopf family of algebras, and while
$\eta \rightarrow 0$ the family structure become trivial. Another
unusual feature of the algebras \Alg is that, under the current
realization, the generating currents corresponding to positive
and negative roots are deformed {\em differently}.
Despite their unusual mathematical features, the algebras \Alg
are believed to have important applications in integrable
quantum field theories such as Sine-Gordon and affine-Toda field
theories. Moreover the study of such kind of algebras would
provide better understanding to the ($(\hbar,\xi)$-) deformed
Virasoro and $W$ algebras which are recently under active study.

In this paper we are motivated to study both the elliptic generalizations
of the algebras \Alg (or the pre-scaling algebras) and their associated infinite Hopf
family of algebras.

The search for elliptic quantum algebras has been lasted for quite some years.
Various elliptic deformed algebras have emerged in several
different contexts, among them there are the Sklyanin algebras of type $sl_N$,
the algebra \Apq \cite{Foda1,Foda2} and its generalization to the $sl_N$ case, \Apqn,
which forms the class of so-called
vertex type elliptic algebras, and the ``elliptic
quantum groups'' of Felder et al \cite{felder1,felder2,felder3}
and the dynamical twisted algebra of Hou et al \cite{hy1} which form the class of so-called
face type elliptic algebras. The above mentioned
elliptic algebras are all realized through the (vertex and face type) Yang-Baxter
relations. The difference between Sklyanin algebra and the algebras \Apqn,
as well as that between Felder et al's algebra and Hou et al's one lie in that,
the modulus for the elliptic entries of $R$-matrices are the same
for the former and are different for the latter algebras. Other examples of
elliptic algebras are the algebras for the deformed screening currents for
the quantum deformed $W$ algebras (defined for any simply-laced underlying
Lie algebra $g$) of the first two of the present authors \cite{qw}
and Konno's algebra $U_{q,p}(\widehat{sl}_2)$ \cite{konno}.

We note that the classification for the
elliptic deformed algebras seems far from complete yet. For example, the last two
types of algebras are realized as current algebras only and their possible Yang-Baxter
type realization are still unknown. Moreover, though the co-structures, or more
explicitly, the quasi-Hopf structures, for the vertex and face type algebras realized
through Yang-Baxter type relations has recently been clarified due to the work of
Fronsdal and Jombo et al, similar structures in the current algebras of \cite{qw}
are still unknown.

In this paper, we shall present a new type of elliptic current algebras
which we denote as \Eg (where $g$ can be any classical simply-laced Lie algebra)
and study the associated infinite Hopf family of algebras structures.
It turns out that the algebras \Eg are quite similar to the algebras of
modified screening currents for the quantum $(q,p)$-deformed $W$-algebras
mentioned above in the level 1 bosonic representations. 
The only difference lies in that, for the algebras
\Eg at level 1, the deformation parameter $q$ is the inverse of the one
in the algebras defined in \cite{qw} (whilst the parameter $\tilde{q}$ is kept
unchanged), and we assume here that $|q|<1$, which
corresponds to $|q|>1$ in \cite{qw} (the algebras in \cite{qw}, however,
were defined only for $|q|<1$ implicitly).
This slight difference prevented us from defining the somewhat well-expected
structure of infinite Hopf family of algebras in \cite{qw}.

The organization of this paper is as follows. In section 2, we shall
give a definition for the current algebra \Eg. Section 3 is devoted to
the study of the structure of associated infinite Hopf family of algebras. In Section 4
we give the bosonic representation for the current algebras \Eg at level 1.
The final section--Section 5--is for some
concluding remarks.

\section{The elliptic current algebra \Eg}

We first give the definition for the elliptic current algebras \Eg.

\begin{defn}
The elliptic current algebra \Eg is the associative algebra generated
by the currents $E_i(z)$, $F_i(z)$, $H^\pm_i(z)$ with $i=1, 2, ..., \mbox{rank} (g)$,
central element $c$ and the unit element 1 with the following relations,

\begin{eqnarray}
H^\pm_i(z)H^\pm_j(w)&=&
\frac{\theta_q\left(\frac{z}{w}p^{A_{ij}/2}\right)
\theta_{\tilde{q}}\left(\frac{z}{w}p^{-A_{ij}/2}\right)}
{\theta_q\left(\frac{z}{w}p^{-A_{ij}/2}\right)
\theta_{\tilde{q}}\left(\frac{z}{w}p^{A_{ij}/2}\right)}
H^\pm_j(w)H^\pm_i(z), \label{39}\\
H^+_i(z)H^-_j(w)&=&
\frac{\theta_q\left(\frac{z}{w}p^{(A_{ij}+c)/2}\right)
\theta_{\tilde{q}}\left(\frac{z}{w}p^{-(A_{ij}+c)/2}\right)}
{\theta_q\left(\frac{z}{w}p^{-(A_{ij}-c)/2}\right)
\theta_{\tilde{q}}\left(\frac{z}{w}p^{(A_{ij}-c)/2}\right)}
H^-_j(w)H^+_i(z), \\
H^+_i(z)E_j(w)&=& (-1)^{A_{ij}} (p^{-{A_{ij}/2}})
\frac{\theta_q\left(\frac{z}{w}p^{A_{ij}/2} p^{c/4}\right)}
{\theta_q\left(\frac{z}{w}p^{-A_{ij}/2} p^{c/4}\right)}
E_j(w)H^+_i(z),\\
H^-_i(z)E_j(w)&=& (-1)^{A_{ij}} (p^{-A_{ij}/2})
\frac{\theta_q\left(\frac{z}{w}p^{A_{ij}/2} p^{-c/4}\right)}
{\theta_q\left(\frac{z}{w}p^{-A_{ij}/2} p^{-c/4}\right)}
E_j(w)H^-_i(z),\\
H^+_i(z)F_j(w)&=& (-1)^{A_{ij}} (p^{A_{ij}/2})
\frac{\theta_{\tilde{q}}\left(\frac{z}{w}p^{-A_{ij}/2} p^{-c/4}\right)}
{\theta_{\tilde{q}}\left(\frac{z}{w}p^{A_{ij}/2} p^{-c/4}\right)}
F_j(w)H^+_i(z),\\
H^-_i(z)F_j(w)&=& (-1)^{A_{ij}} (p^{A_{ij}/2})
\frac{\theta_{\tilde{q}}\left(\frac{z}{w}p^{-A_{ij}/2} p^{c/4}\right)}
{\theta_{\tilde{q}}\left(\frac{z}{w}p^{A_{ij}/2} p^{c/4}\right)}
F_j(w)H^-_i(z),\\
E_i(z)E_j(w)&=& (-1)^{A_{ij}} (p^{-A_{ij}/2})
\frac{\theta_q\left(\frac{z}{w}p^{A_{ij}/2} \right)}
{\theta_q\left(\frac{z}{w}p^{-A_{ij}/2} \right)}
E_j(w)E_i(z),\\
F_i(z)F_j(w)&=& (-1)^{A_{ij}} (p^{A_{ij}/2})
\frac{\theta_{\tilde{q}}\left(\frac{z}{w}p^{-A_{ij}/2} \right)}
{\theta_{\tilde{q}}\left(\frac{z}{w}p^{A_{ij}/2} \right)}
F_j(w)F_i(z),\\
{}[ E_i(z), F_j(w) ] &=&\frac{\delta_{ij}}{(p-1)zw} \left[
\delta\left(\frac{z}{w}q^c\right)H^+_i(zq^{c/2})
-\delta\left(\frac{w}{z}\tilde{q}^{-c}\right)H^-_i(w\tilde{q}^{-c/2}) \right],\\
E_i(z_1)E_i(z_2)E_j(w) \!\!\!\!&-&\!\!\!\! f_{ij}^{(q)}(z_1/w,z_2/w)
E_i(z_1)E_j(w)E_i(z_2) + E_j(w)E_i(z_1)E_i(z_2)\nonumber\\
& & + (\mbox{replacement } ~z_1 \leftrightarrow z_2) =0,
\hspace{0.5cm} A_{ij}=-1,\\
F_i(z_1)F_i(z_2)F_j(w) \!\!\!\!&-&\!\!\!\! f_{ij}^{(\tilde{q})}(z_1/w,z_2/w)
F_i(z_1)F_j(w)F_i(z_2) + F_j(w)F_i(z_1)F_i(z_2) \nonumber\\
& & + (\mbox{replacement } ~z_1 \leftrightarrow z_2) =0,
\hspace{0.5cm} A_{ij}=-1, \label{51}
\end{eqnarray}

\noindent where

\begin{eqnarray*} 
& &f_{ij}^{(a)}(z_1/w,z_2/w) =
\frac{\left(\psi^{(a)}_{ii}\left(\frac{z_2}{z_1}\right)+1 \right)
\left(\psi^{(a)}_{ij}\left(\frac{w}{z_1}\right)
\psi^{(a)}_{ij}\left(\frac{w}{z_2}\right)+1 \right)}
{\psi^{(a)}_{ij}\left(\frac{w}{z_2}\right)
+\psi^{(a)}_{ii}\left(\frac{z_2}{z_1}\right)
\psi^{(a)}_{ij}\left(\frac{w}{z_1}\right)},~~~a=q,~\tilde{q},\\
& &\psi^{(q)}_{ij}(x) =  (-1)^{A_{ij}} p^{-A_{ij}/2}
\frac{\theta_q\left(x^{-1}p^{A_{ij}/2} \right)}
{\theta_q\left(x^{-1}p^{-A_{ij}/2} \right)},\\
& &\psi^{(\tilde{q})}_{ij}(x) =  (-1)^{A_{ij}} p^{A_{ij}/2}
\frac{\theta_{\tilde{q}}\left(x^{-1}p^{-A_{ij}/2} \right)}
{\theta_{\tilde{q}}\left(x^{-1}p^{A_{ij}/2} \right)},
\end{eqnarray*}

\noindent $q$ and $p$ are a pair of deformation parameters with
norms $|q|<1$ and $|p|<1$, $z,w$ are
spectral parameters, $\tilde{q}$ and $q$ are connected by the relation

\begin{eqnarray*}
\tilde{q}/q=p^c,
\end{eqnarray*}

\noindent and $\theta_q(z)$ is the standard elliptic function given by

\begin{eqnarray*}
& & \theta_q(z) =(z|q)_\infty (qz^{-1}|q)_\infty (q| q)_\infty,\\
& &(z| q_1, ..., q_m)_\infty = \prod_{i_1,i_2,...,i_m=0}^\infty
(1-zq_1^{i_1} q_2^{i_2} ... q_m^{i_m}).
\end{eqnarray*}
\end{defn}

Quite analogous to the case of \Alg, the elliptic current algebra given above
enjoys the following features,

\begin{itemize}
\item it has {\em two} deformation parameters $p,q$ and the
``positive'' and ``negative'' currents $E(z)$ and $F(z)$ are deformed differently
(each corresponds to one of the two parameters $q$ and $\tilde{q}$ respectively);
\item the currents $H^\pm_i(z)$ do not commute with themselves in contrast to  the
$q$-affine and Yangian cases.
\end{itemize}

\noindent These features are also shared by the algebras
${\cal A}_{q,p}(\widehat{sl}_N)$, ${\cal A}_{q,p,\hat{\pi}}(\widehat{sl}_2)$
and $U_{q,p}(\widehat{sl}_2)$.

The second feature has a rather significant consequence.
If one consider the subalgebras generated by the currents $H^+(z), E(z)$ or
$H^-(z), F(z)$, it would turn out that they do not form nilpotent
or even solvable subalgebras. However, in the
$q$-affine and Yangian cases similar subalgebras are indeed solvable and with the aid
of a properly defined Manin pairing, they give rise to the structure of quantum
doubles. The non-solvability of such subalgebras in our case might imply that the
algebra \Eg under consideration does not have a simple quantum double structure.

In order to show the more deep relationship between our algebra and the algebras 
\Alg, we give the fillowing proposition which show that the algebra \Eg is an elliptic
extension of \Alg.

\begin{prop}
In the scaling limit

\begin{eqnarray*}
& &p=\mbox{e}^{\epsilon \hbar},~~~q=\mbox{e}^{\frac{\epsilon}{\eta}},~~~
z=\mbox{e}^{i\epsilon u}\\
& & \epsilon \rightarrow 0
\end{eqnarray*}

\noindent the algebra \Eg tends to the algebra \Alg defined in \cite{f}.
\end{prop}

We remark that for the case $g=sl_2$, both the algebra
${\cal A}_{q,p}(\widehat{sl}_2)$ and $U_{q,p}(\widehat{sl}_2)$ would yield
\Al in the scaling limit. Therefore our algebra \Eg has the same scaling
limit as those two algebras for the special underlying Lie algebra $g=sl_2$.
However, for general simply-laced $g$, our algebra \Eg is the only known
algebra which tends to \Alg in the scaling limit. Actually, the generalization
of ${\cal A}_{q,p}(\widehat{sl}_2)$ to the case of $D, E$ series of Lie algebras
are not known to exist. Likewise, the generalization of
$U_{q,p}(\widehat{sl}_2)$ to any other $g$ is also not known to exist. (We noticed the
similarity between our algebra at $g=sl_2$ and $U_{q,p}(\widehat{sl}_2)$. It is
possible that these two algebras are isomorphic, however we do not make this claim
because we did not make it out yet.
\footnote{To compare with \cite{konno}, we one should bare in mind that the the following
change of notations should be made: $q \rightarrow p$, $p \rightarrow q^{2}$ and
$c \rightarrow -c$.})

Another remark is in order here. The algebra \Eg, as well as \Alg defined in \cite{f},
should be regarded as {\em current} algebras only since we do not know the corresponding
Yang-Baxter type realizations. Actually, given a Yang-Baxter type relation one can
define an associative algebra which is certain deformation of the universal enveloping
algebra of some underlying Lie algebra, and,
due to the well-known Ding-Frenkel homomorphism,
one can find a corresponding current realization which is of important usage
for the construction  of infinite dimensional representations. However, the inverse
to Ding-Frenkel homomorphism is some Riemann problem which often does not possesses a
unique solution \cite{khorosh}. Therefore,
given the definition of a current algebra such as
\Eg, one actually cannot associate a unique Yang-Baxter type relation without
putting in extra constraints. It seems quite possible that both the vertex type and
face type elliptic algebras can be obtained from the same current algebra
\Eg by introducing different sets of constraints which lead to different solutions
to the Riemann problem. We hope to consider this problem in later studies.

\section{The structure of infinite Hopf family of algebras for \Eg}

The algebra \Eg defined in the last section
is in fact the representative of an infinite Hopf family
of elliptic algebras which we now specify.

Let $\{ {\cal A}_n,~n\in Z\}$ be a family of associative algebras over $C$ with unit.
Let $\{v^{(n)}_i,~i=1,~...,~\mbox{dim}({\cal A}_n)\}$ be a basis of ${\cal A}_n$.
The maps

\begin{eqnarray*}
\tau_n^{\pm}: {\cal A}_n &\rightarrow& {\cal A}_{n \pm 1}\\
v^{(n)}_i &\mapsto& v^{(n\pm 1)}_i
\end{eqnarray*}

\noindent are morphisms from ${\cal A}_n$ to ${\cal A}_{n \pm 1}$.
For any two integers $n,~m$ with $n<m$, we can specify a pair of morphisms

\begin{eqnarray}
& &\tau^{(m,n)} =Mor({\cal A}_m,~{\cal A}_n) \equiv \tau_{m-1}^+...\tau_{n+1}^+\tau_n^+:~~
{\cal A}_n \rightarrow {\cal A}_m, \nonumber\\
& &\tau^{(n,m)} =Mor({\cal A}_n,~{\cal A}_m) \equiv \tau_{n+1}^-...\tau_{m-1}^-\tau_m^-:~~
{\cal A}_m \rightarrow {\cal A}_n   \label{taumn}
\end{eqnarray}

\noindent with $\tau^{(m,n)}\tau^{(n,m)}=id_{m},~\tau^{(n,m)}\tau^{(m,n)}=id_{n}$.
Clearly the morphisms $\tau^{(m,n)},~n,m\in Z$ satisfy the associativity
condition $\tau^{(m,p)}\tau^{(p,n)} =\tau^{(m,n)}$ and thus make the family of algebras
$\{ {\cal A}_n,~n\in Z\}$ into a category.

\begin{defn}
The category of algebras $\{ {\cal A}_n,~\{\tau^{(n,m)}\},~n,m\in Z\}$ is called an
infinite Hopf family of algebras if on each object ${\cal A}_n$ of the category
one can define the morphisms $\Delta^\pm_n: {\cal A}_n \rightarrow {\cal A}_n
\otimes {\cal A}_{n\pm 1}$, $\epsilon_n: {\cal A}_n \rightarrow C$ and antimorphisms
$S^\pm_n: {\cal A}_n \rightarrow {\cal A}_{n\pm 1}$ such that the following axioms hold,

\begin{itemize}
\item $(\epsilon_n \otimes id_{n+1}) \circ \Delta_n^+ = \tau_n^+,~
(id_{n-1} \otimes  \epsilon_n ) \circ \Delta_n^- = \tau_n^-$ \hfill{(a1)}

\item $m_{n+1} \circ (S^+_n \otimes id_{n+1}) \circ \Delta_n^+
= \epsilon_{n+1} \circ \tau^+_n,~
m_{n-1} \circ (id_{n-1} \otimes S^-_n) \circ \Delta_n^-
= \epsilon_{n-1} \circ \tau^-_n $ \hfill{(a2)}

\item $(\Delta_n^- \otimes id_{n+1}) \circ \Delta_n^+ =
(id_{n-1} \otimes \Delta_n^+ ) \circ \Delta_n^-$ \hfill{(a3)}
\end{itemize}

\noindent in which $m_{n}$ is the algebra multiplication for ${\cal A}_n$.
\end{defn}

\begin{rem}
We remark here that the presentation of infinite Hopf family of algebras
is slightly different from that of
\cite{f} in the trigonometric case. However, the statement that the algebra \Alg
is a representative of an infinite Hopf family of trigonometric algebras still hold true
under the present definition of infinite Hopf family of algebras.
\end{rem}

Let ${\cal A}$ be an associative algebra over $C$ with unit. A trivial example of
infinite Hopf family of algebras is given by the category of algebras $\{ {\cal A}_n \equiv {\cal A},~
\{\tau^{(n,m)}\equiv id_{\cal A}\},~n,m\in Z\}$  with $\Delta_n^\pm$, $\epsilon_n$
and $ S_n^\pm $ identified as the standard Hopf algebra
structures over ${\cal A}$. This trivial
example shows that the infinite Hopf family of algebras can be regarded as some deformation of the standard Hopf
algebra structure. The maps $\Delta_n^\pm$, $\epsilon_n$ and $S_n^\pm$ in the
infinite Hopf family of algebras
are called comultiplications, counits and antipodes by this analogy.

Now let us consider the infinite Hopf family of algebras structure of our
algebra \Eg. For this purpose we
introduce some notations. First, we denote the algebra \Eg by \Egc{}{c}, specifying
explicitly the central extension $c$. We see that this algebra is determined uniquely
as a current algebra by the defining relations (\ref{39}-\ref{51}) provided the following data are
fixed: $g,~q,~p,~c$. In general, given a series of $c_n,~n\in Z$, we can define

\begin{eqnarray*}
q^{(n+1)}/q^{(n)}=p^{c_n},
\end{eqnarray*}

\noindent starting from the data $q^{(1)}=q$, $c_1=c$. It is obvious that
$\tilde{q}=q^{(2)}$ and hence $\tilde{q}^{(n)}=q^{(n+1)}$. We collect the
family of algebras \{\Egc{(n)}{c_n}, $n\in Z$\} where
\Egc{(n)}{c_n} is the algebra \Egc{}{c} with $q$
replaced by $q^{(n)}$ and $c$ by $c_n$. The generating currents
$H_i^\pm(z),~E_i(z)$ and $F_i(z)$ for the algebra \Egc{(n)}{c_n} are denoted as
$H_i^\pm(z; q^{(n)}),~E_i(z; q^{(n)})$ and $F_i(z; q^{(n)})$ etc.

The family of algebras \{\Egc{(n)}{c_n}, $n\in Z$\} can be easily
turned into a category by introducing the morphisms $\tau_n^\pm$

\begin{eqnarray*}
\tau_n^\pm :{\Eegc{(n)}{c_n}} &\rightarrow& {\Eegc{(n\pm 1)}{c_{n\pm 1}}} \\
H^\pm_i(z; q^{(n)}) &\mapsto& H^\pm_i(z; q^{(n\pm 1)})\\
E_i(z; q^{(n)}) &\mapsto& E_i(z; q^{(n\pm 1)})\\
F_i(z; q^{(n)}) &\mapsto& F_i(z; q^{(n\pm 1)})\\
c_n &\mapsto& c_{n\pm 1}
\end{eqnarray*}

\noindent and defining the compositions $\tau^{(n,m)}$ as did in (\ref{taumn}).

The following proposition is one of our major results.

\begin{prop}
The category of algebras $\{\Eegc{(n)}{c_n}, ~\{\tau^{(n,m)}\}, n,m\in Z\}$
form an (elliptic) infinite Hopf family of algebras with the Hopf family structures
given as follows:

\begin{itemize}
\item the comultiplications $\Delta_n^\pm$:
\begin{eqnarray*}
\Delta_n^+ c_n \!\!\!\!&=&\!\!\!\! c_n + c_{n+1}, \label{CoB} \\
\Delta_n^+ H_i^+(z;q^{(n)}) \!\!\!\!&=&\!\!\!\!
H_i^+(z p^{c_{n+1}/4}; q^{(n)})
\otimes H_i^+(z p^{-c_{n}/4}; q^{(n+1)}),\\
\Delta_n^+ H_i^-(z;q^{(n)}) \!\!\!\!&=&\!\!\!\!
H_i^-(z p^{-c_{n+1}/4}; q^{(n)})
\otimes H_i^-(z p^{c_{n}/4}; q^{(n+1)}),\\
\Delta_n^+ E_i(z;q^{(n)}) \!\!\!\!&=&\!\!\!\! E_i(z; q^{(n)})
\otimes 1 + H^-_i(z p^{c_n/4}; q^{(n)}) \otimes
E_i(z p^{c_n/2}; q^{(n+1)}),\\
\Delta_n^+ F_i(z;q^{(n)}) \!\!\!\!&=&\!\!\!\! 1 \otimes
F_i(z; q^{(n+1)}) + F_i(z p^{c_{n+1}/2}; q^{(n)})
\otimes H^+_i(z p^{c_{n+1}/4}; q^{(n+1)}),\\
\\
\\
\Delta_n^- c_n \!\!\!\!&=&\!\!\!\! c_{n-1} + c_n,\\
\Delta_n^- H_i^+(z;q^{(n)}) \!\!\!\!&=&\!\!\!\!
H_i^+(z p^{c_{n}/4}; q^{(n-1)})
\otimes H_i^+(z p^{-c_{n-1}/4}; q^{(n)}),\\
\Delta_n^- H_i^-(z;q^{(n)}) \!\!\!\!&=&\!\!\!\!
H_i^-(z p^{-c_{n}/4}; q^{(n-1)})
\otimes H_i^-(z p^{c_{n-1}/4}; q^{(n)}),\\
\Delta_n^- E_i(z;q^{(n)}) \!\!\!\!&=&\!\!\!\! E_i(z; q^{(n-1)})
\otimes 1 + H^-_i(z p^{c_{n-1}/4}; q^{(n-1)}) \otimes
E_i(z p^{c_{n-1}/2}; q^{(n)}),\\
\Delta_n^+ F_i(z;q^{(n)}) \!\!\!\!&=&\!\!\!\! 1 \otimes
F_i(z; q^{(n)}) + F_i(z p^{c_{n}/2}; q^{(n-1)})
\otimes H^+_i(z p^{c_{n}/4}; q^{(n)});
\end{eqnarray*}

\item the counits $\epsilon_n$:
\begin{eqnarray*}
\epsilon_n ( c_n )\!\!\!\!&=&\!\!\!\!0,\\
\epsilon_n ( 1_n ) \!\!\!\!&=&\!\!\!\! 1,\\
\epsilon_n ( H^\pm_i(z; q^{(n)}))  \!\!\!\!&=&\!\!\!\! 1,\\
\epsilon_n ( E_i(z; q^{(n)})) \!\!\!\! &=& \!\!\!\! 0,\\
\epsilon_n ( F_i(z; q^{(n)})) \!\!\!\! &=& \!\!\!\! 0;
\end{eqnarray*}

\item the antipodes $S_n^\pm$:

\begin{eqnarray*}
S_n^\pm c_n \!\!\!\!&=&\!\!\!\! - c_{n \pm 1},\\
S_n^\pm H^+_i(z; q^{(n)})\!\!\!\!& =&\!\!\!\!
[ H^+_i(z; q^{(n \pm 1)} )]^{-1},\\
S_n^\pm H^-_i(z; q^{(n)})\!\!\!\!& =&\!\!\!\!
[ H^-_i(z; q^{(n \pm 1)} )]^{-1},\\
S^\pm_n E_i(z;q^{(n)}) \!\!\!\!&=&\!\!\!\!
- H^-_i(z p^{-c_{n \pm 1}/4}; q^{(n \pm 1)})^{-1}
E_i(z p^{-c_{n \pm 1}/2}; q^{(n \pm 1)}),\\
S^\pm_n F_i(z;q^{(n)}) \!\!\!\!&=&\!\!\!\!
- F_i(z p^{-c_{n \pm 1}/2}; q^{(n \pm 1)})
H^+_i(z p^{-c_{n \pm 1}/4}; q^{(n \pm 1)})^{-1}.
\end{eqnarray*}
\end{itemize}
\end{prop}

The proof for this proposition is by straightforward calculations.

\begin{rem}
The comultiplications, counits and antipodes given above are analogous to the Drinfeld
Hopf structures for $q$-affine algebras. The difference lies in that, instead of sending
elements of the algebra ${\cal A}_n$ into the tensor product space
of the same algebra, the comultiplications $\Delta^\pm_n$ now send elements
of ${\cal A}_n$ into the tensor product spaces
${\cal A}_n \otimes {\cal A}_{n+1}$ and ${\cal A}_{n-1} \otimes {\cal A}_{n}$
respectively of two neighboring algebras in the family. The shift
in the suffices in the notations of target spaces indicate the crucial difference
between nontrivial infinite Hopf family of algebras and trivial ones.
\end{rem}

In order to understand the meaning of the unusual shift of suffices mentioned above,
we present here another proposition which was first found in \cite{f} in the
trigonometric case.

\begin{prop}
The comultiplication $\Delta^+_n$ induces an algebra homomorphism

\begin{eqnarray*}
\rho: \Eegc{(n)}{c_n+c_{n+1}} &\rightarrow& \Eegc{(n)}{c_n} \otimes
\Eegc{(n+1)}{c_{n+1}}\\
X &\mapsto& \Delta^+_n \tilde{X},
\end{eqnarray*}

\noindent where $X \in \Eegc{(n)}{c_n+c_{n+1}},~~\tilde{X} \in \Eegc{(n)}{c_n}$
and

\begin{eqnarray*}
\tilde{X}=\left\{
\begin{array}{l}
$$c_n$$\\
$$H^\pm_i(z;q^{(n)})$$\\
$$E_i(z;q^{(n)})$$\\
$$F_i(z;q^{(n)})$$
\end{array}
\right.
\quad \quad \quad
\mbox{if}
\quad \quad \quad
X=\left\{
\begin{array}{l}
$$c_n+c_{n+1}$$\\
$$H^\pm_i(z;q^{(n)})$$\\
$$E_i(z;q^{(n)})$$\\
$$F_i(z;q^{(n)})$$
\end{array}
\right. .
\end{eqnarray*}

\noindent Likewise, the comultiplication $\Delta^-_n$ induces an algebra homomorphism

\begin{eqnarray*}
\bar{\rho}: \Eegc{(n-1)}{c_{n-1}+c_{n}} &\rightarrow& \Eegc{(n-1)}{c_{n-1}} \otimes
\Eegc{(n)}{c_{n}}\\
X &\mapsto& \Delta^+_n \tilde{X},
\end{eqnarray*}

\noindent where $X \in \Eegc{(n-1)}{c_{n-1}+c_{n}},~~\tilde{X} \in \Eegc{(n)}{c_n}$
and

\begin{eqnarray*}
\tilde{X}=\left\{
\begin{array}{l}
$$c_{n}$$\\
$$H^\pm_i(z;q^{(n)})$$\\
$$E_i(z;q^{(n)})$$\\
$$F_i(z;q^{(n)})$$
\end{array}
\right.
\quad \quad \quad
\mbox{if}
\quad \quad \quad
X=\left\{
\begin{array}{l}
$$c_{n-1}+c_{n}$$\\
$$H^\pm_i(z;q^{(n-1)})$$\\
$$E_i(z;q^{(n-1)})$$\\
$$F_i(z;q^{(n-1)})$$
\end{array}
\right. .
\end{eqnarray*}

\end{prop}

\begin{cor}
Let $m$ be a positive integer. The iterated comultiplication
$\Delta^{(m)+}_n=(id_n \otimes id_{n+1} \otimes ...
\otimes id_{n+m-2} \otimes \Delta_{n+m-1}^+)
\circ (id_n \otimes id_{n+1} \otimes ...
\otimes id_{n+m-3} \otimes \Delta_{n+m-2}^+) ...
\circ (id_n \otimes \Delta_{n+1}^+) \circ \Delta_n^+$ induces an algebra homomorphism
$\rho^{(m)}$

\begin{eqnarray*}
\rho^{(m)}: \Eegc{(n)}{c_n +c_{n+1} + ... + c_{n+m}} \rightarrow
\Eegc{(n)}{c_n} \otimes \Eegc{(n+1)}{c_{n+1}} \otimes ... \otimes
\Eegc{(n+m)}{c_{n+m}}
\end{eqnarray*}

\noindent in the spirit of Proposition 3.5.
\end{cor}

\begin{rem}
We stress here that the maps $\rho,~\bar{\rho}$ and $\rho^{(m)}$ are {\em algebra
homomorphisms}, whilst $\tau_n^\pm,~\tau^{(n,m)}$ and $\Delta_n^\pm$ etc are
only algebra morphisms. The difference between algebra morphisms and algebra
homomorphisms lies in that the latter preserves the structure functions whilst the
former does not.
\end{rem}

The above proposition and its corollary shows that although the algebras
\Egc{(n)}{c_n} are not co-closed, the tensor product representations can still
be defined using the comultiplications $\Delta^\pm_n$. In particular, if the $c_n=1$
representation for \Egc{(n)}{c_n} is available, then the representations at any
positive integer $c_n$ are all available, left aside the reducibility problems
of such representations. This statement is of particular importance when one
need to realize that the algebra \Egc{}{c} is non-empty for $c \in Z_+\backslash\{1\}$.

\section{Free boson realization of the algebra \Eg at $c=1$}

Having established the infinite Hopf family of algebras structure of
the elliptic current algebra
\Eg, we now turn to consider its simplest infinite dimensional
representation, i.e. the free boson realization at $c=1$.

First we introduce the Heisenberg algebra ${\cal H}_{q,p}(g)$ with generators
$a_i[n],~P_i,~Q_i, i=1,...,$ $\mbox{rank}(g), n\in Z\backslash \{0\}$ and
generating relations

\begin{eqnarray*}
& &[ a_i[n],~a_j[m] ] = \frac{1}{n}
\frac{(1-q^{-n})(p^{nA_{ij}/2}-p^{-nA_{ij}/2})(1-(pq)^n)}{1-p^n} \delta_{n,m},\\
& &[ P_i,~Q_j ] = A_{ij},
\end{eqnarray*}

\noindent where $A_{ij}$ is the Cartan matrix for the Lie algebra $g$.
Let

\begin{eqnarray*}
s_i^+[n]=\frac{a_i[n]}{q^n-1},~~~~s_i^-[n]=-\frac{a_i[n]}{(pq)^{-n}-1}
\end{eqnarray*}

\noindent and define the (deformed) free boson fields

\begin{eqnarray*}
\varphi_i(z)=\sum_{n \neq 0} s_i^+[n] z^{-n},~~~~
\psi_i(z) = \sum_{n \neq 0} s_i^-[n] z^{-n}.
\end{eqnarray*}

\begin{prop}
The following bosonic expressions give a level $c=1$ realization for the algebra
\Eg on the Fock space of the Heisenberg algebra ${\cal H}_{q,p}(g)$,

\begin{eqnarray*}
& & E_i(z)=e^{Q_i} z^{P_i} : \exp[ \varphi_i(z(pq)^{1/2})] :,\\
& & F_i(z)=e^{-Q_i} z^{-P_i} : \exp[ -\psi_i(zq^{1/2})] :,\\
& & H^+_i(z)= :E_i(zp^{1/4}) F_i(zp^{-1/4}):,\\
& & H^-_i(z) = :E_i(zp^{-1/4}) F_i(zp^{1/4}):,
\end{eqnarray*}

\noindent where : : means taking all subexpressions consisted of $a_i[n]$ with $n>0$ and
$P_i$ to the right of expressions consisted of $a_i[n]$ with $n<0$ and
$Q_i$.
\end{prop}

The proof for this proposition is also by straightforward but tedious calculations.

\section{Concluding Remarks}

In this paper we obtained the new elliptic current algebras \Eg and
showed that these algebras has a structure of infinite Hopf family of
algebras.
So far we have obtained two kinds of nontrivial infinite Hopf family of algebras:
trigonometric (for the algebras
\Alg) and elliptic (for the algebras \Eg). It is thus an interesting question to
ask whether there exists any rational algebras which has the same co-algebraic
structure.

It is interesting to mention that the comultiplications appearing in
such co-structures are all of the Drinfeld type, which closes over the currents
themselves and does not require the resolution to the inverse problem (Riemann problem)
of the Ding-Frenkel homomorphism.  Recall that two kinds of comultiplications (and thus
two kinds of Hopf algebra structures) are known for the standard $q$ affine algebras.
The algebras \Alg and \Eg should be considered as some deformation of $q$-affine algebras
and under such deformations the difference between the two Hopf algebra structures for
$q$-affine algebras become clear: the standard Hopf structure for $q$-affine algebras
is inherited into the Yang-Baxter type realizations for the algebras
\Apq and ${\cal B}_{q,\lambda}(\hat{g})$ \cite{jimbonew}
and define the quasi-triangular quasi-Hopf structures in those algebras, and the
Drinfeld type Hopf structure is inherited into the current realizations for the algebras
\Alg and \Eg and gives rise to the structure of infinite Hopf family of algebras.
The relation between the quasi-triangular quasi-Hopf structure and infinite Hopf
family of algebras is an interesting open problem to be
answered in later studies.

We should emphasis that this work is only a preliminary study for the
algebras \Eg themselves. Besides the definition and level 1 bosonic realization,
we know very little about these algebras, especially their detailed representation
theory, vertex operators, Yang-Baxter type realizations etc. The physical applications
should also be considered.

Finally, the structures of infinite Hopf family of algebras is still poorly
understood yet.  We do not know whether there exists a quantum double construction over
the infinite Hopf family of algebras and, if not, what kind of new structure will take
the place of the standard quantum doubles. Also, the classical counterpart of the
infinite Hopf family of algebras is unknown and it seems that all these problems
deserve further investigations.


\begin{thebibliography}{40}
\bibitem{felder1}
Felder, G., Conformal field theory and integrable systems associated to elliptic curves,
{\it Proc. ICM Z\"urich 1994, 1247, Birkh\"auser (1994)};
Elliptic quantum groups, {\it Proc. ICMP Paris 1994, 2118, International Press (1995)}.

\bibitem{felder2}
Felder, G., Varchenko, A., On representation theory of the elliptic quantum group
$E_{\\tau,\eta}(sl_2)$, {\tt q-alg/9601003}, {\it Commun. Math. Phys. 181} (1996) 741.

\bibitem{felder3}
Enriquez, B., Felder, G., Elliptic quantum groups $E_{\tau,\eta}(sl_2)$ and
quasi-Hopf algebras, {\tt q-alg/9703018}.


\bibitem{Foda1}
Foda, O., Iohara, K., Jimbo, M., Kedem, R., Miwa, T., Yan, H.,
An elliptic quantum algebra for $\widehat{sl}_2$.
{\it Lett. Math. Phys. 32} (1994) 259--268.

\bibitem{Foda2}
Foda,O., Iohara,K., Jimbo,M., Kedem,R., Miwa,T., Yan, H.,
Notes on highest weight modules of the elliptic algebra \Apq.
{\it Prog. Theoret. Phys., Supplement, 118} (1995) 1--34.

\bibitem{hy1}
Hou, B.-Y., Yang, W.-L., Dynamically twisted algebra $A_{q,p,\hat{\pi}}(\widehat{gl}_2)$
as current algebra generalizing screening currents of $q$-deformed Virasoro algebra,
Preprint {\tt q-alg/9709024}.

\bibitem{f} Hou, B.-Y., Zhao, L., Ding, X.-M., The algebra ${\cal A}_{\hbar,\eta}
(\hat{g})$ and Hopf family of algebras, Preprint {\tt q-alg/9703046},
{\it J. Geom. Phys. (in press)}.

\bibitem{jimbonew}
Jimbo, M., Konno, H., Odake, S., Shiraishi, J., Quasi-Hopf twistors for
elliptic quantum groups, {\tt q-alg/9712029}.

\bibitem{khorosh}
Khoroshkin, S., Lebedev, D., Pakuliak, S., Yangian algebras and classical Riemann
problems, Preprint {\tt ITEP-TH-66/97}.

\bibitem{KLP}
Khoroshkin, S., Lebedev, D., Pakuliak, S., Elliptic algebra \Apq in the
scaling limit, Preprint {\tt q-alg/9702002}.

\bibitem{konno}
Konno, H., An elliptic algebra $U_{q,p}(\widehat{sl}_2)$ and the fusion RSOS model,
Preprint {\tt q-alg/9709013}.

\bibitem{qw} Zhao, L., B.-Y. Hou, Note on the algebra of screening currents for the
quantum deformed $W$ algebras, {\it J. Phys. A: Math. Gen. 30} (1997) 7659.
\end{thebibliography}
\end{document}